\documentclass[12pt,a4paper]{article}
\usepackage[cp1250]{inputenc}
\usepackage[OT1]{fontenc}
\usepackage{amsmath,amssymb}
\usepackage{amsthm}
\usepackage[pdftex]{graphicx}
\usepackage[all]{xy}

\newtheorem{thm}{Theorem}
\newtheorem{lemma}{Lemma}
\newtheorem{df}{Definition}
\newtheorem{remark}{Remark}

\newcommand{\lag}{\mathfrak g}
\newcommand{\lap}{\mathfrak p}
\newcommand{\so}{\mathfrak{so}}
\newcommand{\lah}{\mathfrak h}
\newcommand{\gl}{\mathfrak{gl}}
\newcommand{\R}{\mathbb R}
\newcommand{\N}{\mathbb N}
\newcommand{\C}{\mathbb C}
\newcommand{\V}{\textrm V}
\newcommand{\U}{\textrm U}
\newcommand{\mM}{\textrm M}
\newcommand{\W}{\textrm W}
\newcommand{\E}{\textrm E}
\newcommand{\SO}{\textrm{SO}}
\newcommand{\SL}{\textrm{SL}}

\newcommand{\Spin}{\textrm{Spin}}
\newcommand{\G}{\textrm G}
\newcommand{\LGP}{\textrm P}
\newcommand{\DO}{\mathbf D}
\newcommand{\Sp}{\mathbb S}
\newcommand{\T}{\mathbb T}
\newcommand{\Z}{\mathbb Z}
\newcommand{\PG}{\mathcal G}
\newcommand{\M}{\mathcal M}

\newcommand{\ATB}{\mathcal A}
\newcommand{\dis}{\mathcal H}

\newcommand{\osu}{\mathcal U}

\newcommand{\p}{\partial}

\newcommand{\qt}{\frac{1}{4}}

\newcommand{\ra}{\rightarrow}
\newcommand{\eps}{\epsilon}
\newcommand{\e}{e}

\title{k-Dirac operator and parabolic geometries}
\author{
Tomáš Salač\footnote{Research supported by GACR 201/09/H012 and
SVV-2011-263317.}\\
Charles University, Prague\\
\texttt{salac@karlin.mff.cuni.cz}}

\begin{document}
\maketitle

\abstract{ The principal group of a Klein geometry has
canonical left action on the homogeneous space of the geometry and
this action induces action on the spaces of sections of vector
bundles over the homogeneous space. This paper is about construction
of differential operators invariant with respect to the induced
action of the principal group of a particular type of parabolic
geometry. These operators form sequences which are related to the
minimal resolutions of the $k$-Dirac operators studied in Clifford
analysis.}
\section{Introduction.}
Let $\R_n$ be a Clifford algebra of $\R^n$ with an Euclidean
scalar product and let $\{\varepsilon_j,1\le j\le n\}$ be the standard basis. The $k$-Dirac operator
$\{\partial_1,\ldots,\partial_k\}$ is an over-determined system of
first order differential operators. Let $f$ be a smooth $\R_n$-valued function on $\R^{kn}$. Then
\begin{equation}\label{k-dirac operator}
\partial_{i}f=\sum_{1\le j\le n}\varepsilon_j\partial_{ij}f,
\end{equation}
where we identify $\R^{kn}$ with the space of the real matrices $M(n,k,\R)$ of rank $n\times k$. Then $\partial_{ij}$ are
the usual partial derivatives and $\varepsilon_j$ stands for the multiplication by the Clifford number $\varepsilon_j$.

The solutions of the $k$-Dirac equation $(\forall i:\partial_if=0)$ are called
monogenic functions in several Clifford variables. Monogenic
functions share analogous properties as holomorphic functions in one
complex variable and from this point of view, the $k$-Dirac operator
can be viewed as a generalization of the Cauchy-Riemann operator.

Interesting behaviour of holomorphic functions in $n$ variables on domains in $\C^n$, such as Hartog's paradox, can be characterized by the sheaf cohomology of holomorphic functions. The sheaf cohomology can be defined as the left derived functor to the functor of global sections. Thus one needs a suitable resolution of the sheaf of holomorphic functions. This is usually the Dolbeault resolution. In light of these facts, natural question is to find a resolution of the sheaf of monogenic functions.

Let us first mention some results from Clifford analysis.
From a general theorem for differential operators with constant
coefficients follows existence of a finite resolution of the sheaf monogenic
functions, see \cite{CSSS}. However there is no direct way how to
translate this general theorem into explicit form of operators in the resolution. Candidates
for the minimal resolutions were computed in the papers
\cite{SSS},\cite{SSSL}. If $n\ge 2k$, all the operators in the sequence are polynomial combinations of the operators (\ref{k-dirac operator}) and the formulas do not explicitly depend on $n$. Thus we suppress the parameter $n$ and talk about the $k$-Dirac
operator. This case is called the stable case. The unstable range, i.e. $n<2k$, is more complicated since
exceptional syzygies, which do not come from the commutation relations of the operators (\ref{k-dirac operator}), arise. Methods used in these papers include mostly methods from partial differential equations, homological and computer algebra. These methods run very quickly into very high computational complexity.

For $n=4$ the operator is called the
$k$-Cauchy-Fueter operator. This operator was studied in the papers
\cite{BAS} and \cite{BS} from more geometric point of view, that is
from the point of view of symmetries of the operator using methods of
representation theory. The k-Cauchy-Fueter operator was studied in more geometric setting of quaternionic manifolds in the paper
\cite{B}. In this paper we will study k-Dirac operator in similar geometric setting as in the paper \cite{B} so let us briefly recall some basic facts.

A quaternionic manifold $M$ is a manifold of dimension $4n$ with a
reduction of structure group of the tangent bundle to $Sp(1)Sp(n)$ admitting a compatible torsion
free connection. There is the unique normal parabolic geometry
$(\mathcal G,\omega)$ induced by the quaternionic structure, see
\cite{CS}. To the principle fibre bundle $\mathcal P$ one can
associate a triple of spaces. The Penrose transform transfer
cohomological data from the first space, so called twistor space, over the second space to a sequence of differential operators on the third space, i.e. on the quaterninonic manifold $M$, and to
solutions of the differential operators in the sequence. Vanishing of
the cohomology groups on the twistor space implies local exactness of the sequence on the manifold $M$. In this way, one gets a resolution of the k-Cauchy operator on quaternionic manifolds. For more discussion on comparison between k-Cauchy-Fueter operators living in
quaternionic manifolds and Euclidean spaces see \cite{CSS}. For more information about the Penrose transform see \cite{BE} and
\cite{WW}.

Each operator in the resolutions of k-Cauchy-Fueter operators living in quaternionic manifolds is invariant with respect to the induced parabolic structure $(\mathcal G,\omega)$.
The parabolic structure puts severe conditions on invariant
operators and thus we get rid of ambiguity which was present on the
Euclidean space. Linear differential operators invariant with
respect to parabolic structures has been studied intensively through
the last century. The operators of the first order were completely
classified and characterized in the paper \cite{SS}. After reducing the structure group of a parabolic geometry to its reductive part, any such
operator is given by the derivation with the chosen Weyl connection
and a linear projection. Invariance of such operator gives equation
for the generalized conformal weight and if the equation is
satisfied then the resulting formula does not depend on the chosen
Weyl connection. However situation with higher order operators is
far more complicated since there are usually more operators which
can be combined and finding the right combination is a non-trivial
task.

The starting point for understanding differential operators
invariant with respect to some particular parabolic structure is
understanding $\G$-invariant operators on the homogeneous (flat) space of
the geometry. Let $\G$ be a principal group  of a parabolic
geometry, i.e. $\G$ is a semi-simple Lie group with a parabolic
subgroup $\LGP$. Let $(\V,\rho)$ be a representation of the group
$\LGP$ and let $\G\times_\LGP\V$ be the associated vector bundle.
The space of sections $\Gamma(\G\times_\LGP\V)$ can be identified
with the space $\mathcal C^\infty(\G,\V)^\LGP$ of $\V$-valued
$\LGP$-equivariant functions on the total space $\G$. The action of
the group $\G$ on the space of smooth sections is, using the
isomorphism, defined by
\begin{equation}\label{action on sections}
(g.f)(h)=f(g^{-1}h)
\end{equation}
where $g,h\in\G,f\in\mathcal C^\infty(\G,\V)^\LGP$.
An operator
$$D:\Gamma(\G\times_\LGP\V)\ra\Gamma(\G\times_\LGP\W)$$
is called $\G$-invariant if
\begin{equation}\label{g invariance}
D(g.s)=g.(Ds).
\end{equation}
for any $s\in\Gamma(\G\times_\LGP\V)$ and for all $g\in\G$. It is
well known fact that there is a one-to-one correspondence between
$\G$-invariant operators and homomorphisms of generalized Verma
modules.

In the regular character, there is complete classification of the
homomorphisms and one gets so called BGG-sequences, see \cite{CSlS}. Namely, the
highest weights of generalized Verma modules are connected by the
affine action of the Weyl group of the Lie algebra $\lag$ of the Lie
group $\G$ and comparable with respect to the partial order iff there is a non-zero homomorphism between generalized Verma
modules such that the domain of the homomorphism is the Verma module whose highest weight is smaller or equal to the highest weight of the latter Verma module.

In the singular character holds only the if part.
Nevertheless we have a necessary conditions on pairs of generalized Verma
modules which admit a non-trivial homomorphism. To keep track on sequences of homomorphisms of generalized Verma modules in the singular character, we usually draw so called singular Hasse graphs, i.e. oriented graphs where the vertices are the highest weights of the generalized Verma modules and the arrows correspond to non-trivial homomorphism of generalized Verma modules which may or may not exist. In the picture of invariant differential operator, we know for which vector bundles we should look for.

In \cite{F} were computed singular Hasse graphs which coincide with
the graphs of minimal resolutions which were found by people in Clifford analysis. The singular Hasse graphs correspond to sequences
of first and second order operators. Existence of the first order
operators was proved but existence of the second order operators was in general
only conjectured. Also local form of the operators was
given only in some particular cases. In contrast to the $k$-Cauchy-Fueter operator, the geometry is no
longer 1-graded but is 2-graded which brings some new difficulties.

The main goal of this paper is to fill some gaps into comparison of sequences of operators studied in \cite{F} and the minimal resolutions found by people from Clifford analysis. The main result is construction of $\G$-invariant
second order operators conjectured in \cite{F} in the stable range
on the homogeneous space of the geometry.

\begin{remark}\label{uniqueness of G-inv op}
In general, a $\G$-invariant operator going between two fixed homogeneous vector bundles in one direction need not be unique up to a scalar multiple. However in this case $k=2,3$ uniqueness was proved in the paper \cite{BC}. See also \cite{H}.
\end{remark}

I have used a construction called splitting operators. The splitting operators are given by polynomials in Curved
Casimir operators. The Curved Casimir operator was introduced in \cite{CS}. A formula (\ref{final diff operator}) gives a general second order operator as a linear combination of operators invariant with respect to the Levi factor of the parabolic subgroup. The coefficients are given in the theorem \ref{wght 1}. This formula is given in a preferred Weyl structure and Cartan gauge over an open affine subset of the homogeneous space. An explicit formula is given in the case $k=2$ in a local chart in (\ref{first operator}). I have introduced some additional assumptions on the sections over the affine subset of the homogeneous space, this is related to non-trivial grading of the geometry. In this way, we can consider sections which have the same freedom as in the Euclidean setting. These additional assumptions are used in the last section in the case $k=3$ where one can again verify directly that the sequence is a complex. In both cases we have that the sequences of symbols are exact when restricted to the distribution which rules the geometry. These are the theorems \ref{complex for k=2}. and \ref{complex for k=3}.

The Penrose transform is being used also by L. Krump in \cite{K} in
order to obtain more information about sequences which correspond to
$k$-Dirac resolutions, in particular in the non-stable range. The
Penrose transform might be hopefully the right tool to answer the
question about local exactness of the operators given in the paper

We close the introductory section with structure of the paper. First is introduced the parabolic geometry where the operators
live together with some necessary notation. Then we recall fundamental results from \cite{F} which we will need. Then is introduced the homogeneous vector bundle where all computations will be carries out with small hint how this bundle can be found out. Then we recall the definition of the Curved Casimir operator and give local formulas for the operator using a local adapted frame. We briefly recall how to get splitting operators from the Curved Casimir operator. Then we set preferred choices over the affine open set of the homogeneous space, in particular we will introduce some particular local adapted frame which turns out to be most useful for our purposes. Then we carry out computations in the preferred trivializations to obtain the formula for general second order operator and verify that this construction gives the operators conjectured in \cite{F}. In the last two section, the general formula is given in to local form for $k=2$ and $k=3$. In the case $k=3$, we will work only with the smaller set as mentioned in the previous paragraph.

\section{Parabolic geometry behind the $k$-Dirac operator.}
Let $\G$ be a $4:1$ cover of the connected component of the identity of the real group
$\SO(n+k,k)$\footnote{If $k>2$, the fundamental group of $\SO(n+k,k)$ is
isomorphic to $\Z_2\times\Z_2$ and thus the Lie group $\G$
is the universal cover of $\SO(n+k,k)$. We assume throughout the paper that $k>2$ and $n\ge 2k$. The
fundamental group of $\SO(n+2,2)$ is isomorphic to $\Z\times\Z_2$ and
the group $\G$ is then the covering corresponding to the subgroup
$2\Z\times1\subset\Z\times\Z_2$.}. The group $\G$ has a natural
transitive action on the Grassmannian manifold $V_0(k,n+2k)$ of
oriented maximal isotropic vector subspaces of dimension $k$ in the
vector space $\R^{n+2k}$ with a $\G$-invariant quadratic form of
signature $(n+k,k)$. A parabolic subgroup $\LGP$ is the stabilizer
of a chosen maximal isotropic $k$-dimensional subspace \footnote{The matrix of the scalar product with respect to the preferred basis is 
$$\left(
\begin{array}{ccc}
0&0&1_k\\
0&1_n&0\\
1_k&0&0\\
\end{array}
\right),$$ where block matrices on the diagonal are square matrices
with ranks equal to $k,n,k$. Thus we get a preferred
matrix realization of $\so(n+k,k)$. The preferred
maximal isotropic subspace $L$ is the one spanned by null vectors
$e_i,i=1,\ldots,k$ of the preferred basis. A natural candidate for a complement of $L$ in $L^\bot$ is then the subspace spanned by the vectors $e_{k+1},\ldots,e_{k+n}$ of the preferred basis.} $L$ in $\R^{n+2k}$. Thus we can think of the
homogeneous space $\G/\LGP$ also as the Grassmannian $V_0(k,n+2k)$. A choice of a complement of $L$ in $L^\bot$ gives a natural gradation
\begin{equation}\label{gradation on lie alg}
\mathfrak{g}=\mathfrak{g}_{-2}\oplus\mathfrak{g}_{-1}\oplus\mathfrak{g}_{0}\oplus\mathfrak{g}_{1}\oplus\mathfrak{g}_{2},
\end{equation}
on the Lie algebra $\lag$ of $\G$. 

Let $\G_0$ be the subgroup of $\LGP$ whose adjoint action preserve the gradation (\ref{gradation on lie
alg}). The group $\G_0$ is a maximal reductive subgroup of $\LGP$ and its Lie algebra is
isomorphic to $\gl(k,\R)\oplus\so(n)$. Let us denote by
$\lap_+=\lag_1\oplus\lag_2,\lag_-=\lag_{-2}\oplus\lag_{-1},\lag^i=\oplus_{j\ge
i}\ \lag_j$. We denote the one dimensional center of $\lag_0$ by $\E_0$.

There are isomorphisms of $\G_0$-modules
\begin{equation}\label{g0 pieces in lap}
\mathfrak{g}_{1}\cong\V\otimes\E,\ \mathfrak{g}_{2}\cong\Lambda^{2}\V\otimes\R,
\end{equation}
where $\V$, resp. $\E$, resp. $\R$, denotes the defining representation of
$\gl(k,\R)$, resp. of $\so(n)$, resp. the trivial representation of
$\so(n)$. 

The Lie bracket $\lag_1\otimes\lag_1\ra\lag_2$ is the tensor product of anti-symmetrization and contraction with respect to the $\so(n)$-invariant Euclidean scalar product $g$ on $\E$, i.e.
\begin{eqnarray}\label{lie bracket}
&\wedge\otimes Tr:(\V\otimes\E)\otimes(\V\otimes\E)\ra\Lambda^2\V\otimes\R&\\
&(v_1\otimes\varepsilon_\alpha)\otimes(v_2\otimes\varepsilon_\beta)\mapsto v_1\wedge v_2\otimes g(\varepsilon_\alpha,\varepsilon_\beta)\nonumber
\end{eqnarray}
with $v_1,v_2\in\V,\varepsilon_\alpha,\varepsilon_\beta\in\E$. 

\subsection{Notation for $\lag_0$-modules.}
Let $\lah=\{H\ |\ H=(h_{ij})$ is a trace free diagonal $(k\times
k)$-matrix$\}$ be a Cartan subalgebra of the algebra
$\mathfrak{sl}(k,\R)$. Let $\epsilon_i,1\le i\le k,$ be the element
of the dual space $\lah^\ast$ defined by $\epsilon_i(H)= h_{ii}$. The set
$\{\epsilon_i\}_{i=1}^k$ is a spanning set of the space $\lah^\ast$. Let
$\lambda\in\lah^\ast$, then we can express
$\lambda=\sum_{i=1}^k\lambda_i\epsilon_i$ with $\lambda_i\in\R$. Then we write $\lambda=(\lambda_{1},\ldots,\lambda_{k})$ and we denote by $\lambda_{(ij)}$ the weight $\lambda_{(ij)}=(\lambda_{1},\ldots,\lambda_{i}+1,\ldots,\lambda_{j}+1,\ldots,\lambda_{k})$ with $1\le i<j\le k$.

The weight $\lambda\in\lah^\ast$ is dominant integral iff $\lambda_i-\lambda_j\in\N\cup\{0\}$ for all $1\le i<j\le k$. Let $\V_\lambda$ be an irreducible complex finite dimensional representation of $\mathfrak{sl}(k,\R)$ with a highest weight $\lambda$. We extend the action of $\mathfrak{sl}(k,\R)$ to $\gl(k,\R)$ on $\V_\lambda$ by $1_k\in\gl(k,\R)\mapsto(\sum\limits_i\lambda_i) Id_{\V_\lambda}\in End(\V_\lambda)$.
\begin{remark}
All considered representations are complex, if necessary we consider
complexification of a real representation. Tensor products are over
the field of complex numbers.
\end{remark}

\subsection{Sequences of differential operators related to minimal resolutions of the $k$-Dirac operator.}
Let us fix $n,k$ such that $n\ge2k\ge4$ and let us suppose for simplicity that $n$ is even. The case $n$ is odd
proceeds similarly with obvious modifications. We want to find a second order $\G$-invariant differential operator $D$ which belongs to the sequence starting with the $k$-Dirac operator. From \cite{F} we know that
\begin{equation}\label{diff op}
D:\Gamma(\G\times_{\LGP}\W(\lambda))\ra\Gamma(\G\times_{\LGP}\W(\nu))
\end{equation}
such that highest weights of the $\LGP$-modules dual to the irreducible
modules $\W(\lambda),\W(\nu)$ lie either on the affine orbit of the
weight $\frac{1}{2}(-n+1,\ldots,-n+1|1,\ldots,1)$ or on the affine
orbit of the weight $\frac{1}{2}(-n+1,\ldots,-n+1|1,\ldots,1,-1)$.\footnote{Here $|$ separates the weight of $\mathfrak{gl}(k,\R)$ from
the weight of $\mathfrak{so}(n)$.} This implies that the module $\W(\lambda)$, resp. $\W(\mu)$ is isomorphic to the module
$\V_\lambda\otimes\Sp_\pm,$ resp. $\V_\mu\otimes\Sp_\pm$ where:
\begin{enumerate}
\item $\V_\lambda$, resp. $\V_\mu$ is the irreducible $\mathfrak{gl}(k,\R)$-module with the highest weight $\lambda$, resp. $\mu$.
\item $\Sp_\pm$ is isomorphic to the
complex $\so(n)$-module $\Sp_+$ or to the module $\Sp_-$.
\item $\nu=\lambda_{(ij)}$ for some $1\le i<j\le k$.
\item The weights $\lambda-\frac{1}{2}(n-1,\ldots,n-1)$ and $\lambda_{(ij)}-\frac{1}{2}(n-1,\ldots,n-1)$ have the Young diagrams symmetric with respect to the reflection along the main diagonal.
\end{enumerate}

\subsection{Receipt for finding the operator (\ref{diff op}).}
In order to find the operator (\ref{diff op}), we will use the Curved Casimir operator on a suitable homogeneous vector bundle. First we need to find the right homogeneous vector bundle. Let us consider the following facts.

Let us choose a Weyl structure. With this choice, we can write the operator (\ref{diff op}) as a combination of $\G_0$-invariant operators. It is reasonable to expect that the highest (second) order part of the operator is a combination of
operators which are given by differentiating twice with vector fields lying only in
the distribution $\G\times_\LGP(\lag^{-1}/\lap)$ of the tangent
bundle $\mathrm T\G/\LGP\cong\G\times_\LGP\lag/\lap$ and algebraic $\G_0$-equivariant projections. Thus we have to consider maps from the space of the sections of the homogeneous vector bundle associated to $\V_\lambda\otimes\Sp_\pm$ to the space of the sections of the bundle associated to $\lag_1\otimes\lag_1\otimes\V_{\lambda}\otimes\Sp_\pm$.

We notice that the multiplicity of the target
module $\V_{\lambda_{(ij)}}\otimes\Sp_\pm$ in the $\G_0$-module
$\lag_1\otimes\lag_1\otimes\V_{\lambda}\otimes\Sp_\pm$ is equal to four
and thus there are four second order $\G_0$-invariant differential
operators to combine. But there is, up to a scalar multiple\footnote{See the remark \ref{uniqueness of G-inv op}.}, a unique
combination which can be the highest order part of a $\G$-invariant
operator and the first aim is to find such linear
combination.

Thus we need to find some $\LGP$-module, which is completely reducible as
$\G_0$-module, containing the module $\V_\lambda\otimes\Sp_\pm$ and the target module
$\V_{\lambda_{(ij)}}\otimes\Sp_\pm$, both with the multiplicity one. A natural choice is the minimal
$\LGP$-module $\mM_\lambda^\pm$ in the module
$\mathfrak{g}\otimes\V_{\lambda}\otimes\mathbb{S}_{\pm}$ containing $\G_0$-submodule
$\E_0\otimes\V_{\lambda}\otimes\mathbb{S}_{\pm}$ where $\E_0$ is the
one-dimensional centre of $\lag_0$. It is easy to see that the
$\LGP$-module $\mM^\pm_\lambda$ is, as the $\G_0$-module, isomorphic
to
\begin{equation}\label{decomposition of big module}
\mM^\pm_\lambda\cong_{\G_0}\E_0\otimes(\V_{\lambda}\otimes\mathbb{S}_\pm)\
\oplus\ \lag_{1}\otimes(\V_{\lambda}\otimes\mathbb{S}_\pm)\ \oplus\
\lag_{2}\otimes(\V_{\lambda}\otimes\mathbb{S}_\pm).
\end{equation}
So we see that $\mM^\pm_\lambda$ contains $\V_\lambda\otimes\Sp_\pm$ with multiplicity one. Now we prove that $\mM^\pm_\lambda$ contains $\V_{\lambda_{(ij)}}\otimes\Sp_\pm$ with multiplicity one.
\begin{lemma}
There is a unique $\G_0$-submodule in $\mM^\pm_\lambda$ isomorphic
to $\V_{\lambda_{(ij)}}\otimes\Sp_\pm$ and this module is a
submodule of $\lag_2\otimes\V_\lambda\otimes\Sp_\pm$.
\end{lemma}
Proof: It is easy to see that a submodule of $\mM^\pm_\lambda$
isomorphic to $\V_{\lambda_{(ij)}}\otimes\Sp_\pm$ must be a
submodule of $\lag_2\otimes\V_\lambda\otimes\Sp_\pm$. First let us
notice that
$\lag_2\otimes\V_\lambda\otimes\Sp_\pm\cong\Lambda^2\V\otimes\C\otimes\V_\lambda\otimes\Sp_\pm\cong\Lambda^2\V\otimes\V_\lambda\otimes\Sp_\pm$.
Thus we need to show that $\Lambda^2\V\otimes\V_\lambda$ contains
the representation $\V_{\lambda_{(ij)}}$ with multiplicity one.
\begin{lemma}\label{multiplicity one}
Let $\lambda=(\lambda_{1},\ldots,\lambda_{k})$ and
$\lambda_{(ij)}=(\lambda_{1},\ldots,\lambda_{i}+1,\ldots,\lambda_{j}+1,\ldots,\lambda_{k})$
be integral dominant weights for $\mathfrak{sl}(k,\R)$. Then the
module $\V_{\lambda_{(ij)}}$ appears with multiplicity one in the
module $\Lambda^{2}\V\otimes\V_{\lambda}$.
\end{lemma}
The lemma\ \ref{multiplicity one}. can proved by the Klimyk lemma.
First we need following notation. Let $\lag'$ be a semi-simple Lie
algebra, let $\V'$ be a $\lag'$-module and $\sigma$ be an element of
the weight lattice of $\lag'$. Let $m^{\V'}_{\sigma}$ be the
dimension of the weight space of the weight $\sigma$ in the module
$\V'$. For a regular weight $\sigma$, let $\eta_\sigma$ be the
determinant of $S$ where $S$ is the unique element of the Weyl group
such that $S(\sigma)$ is a dominant weight $[\sigma]$ in the orbit
of $\sigma$ under the action of the Weyl group and set
$\eta_\sigma=0$ for a singular weight.
\begin{lemma}[Klimyk lemma]\label{klimyk lemma}
Let $\lag'$ be a semi-simple Lie algebra. Let $\nu,\nu',\nu''$ be
integral dominant weights for $\lag'$, let $\rho$ be the lowest
weight of $\lag'$. Then the multiplicity $n_{\nu}$ of the
irreducible representation $\V_{\nu}$ with the highest weight $\nu$
in $\V_{\nu'}\otimes\V_{\nu''}$ is equal to $\sum
m^{\V_{\nu'}}_{\sigma}\eta_{\sigma+\nu''+\rho}$ where the sum is
taken over all the weights $\sigma$ of the representation
$\V_{\nu'}$ for which $[\sigma+\nu''+\rho]=\nu+\rho$.
\end{lemma}
Proof: See \cite{H}.

Proof of the lemma \ref{multiplicity one}. In the formulation of lemma
\ref{klimyk lemma}, we have $\nu'=(1,1,0,\ldots),\ \ \nu''=\lambda,
\ \nu=\lambda_{(ij)},\ \ \rho=(k-1,k-2,\ldots,1,0)$. The weights of
$\Lambda^{2}\V$ are of the form $\{\sigma_{\alpha\beta}|1\leq
\alpha<\beta\leq k\}$ where $\sigma_{\alpha\beta}$ has $1$ at
$\alpha$-th and $\beta$-th entry and $0$ otherwise. Moreover
$\lambda+\rho$ and $\lambda_{(ij)}+\rho$ are strictly dominant and
for any $\alpha,\beta$ the weight
$\rho+\lambda+\sigma_{\alpha\beta}$ is dominant (not necessarily
strictly). Since the action of the Weyl group of
$\mathfrak{sl}(k,\R)$ only permutes the entries of a weight, the
only dominant weight in the orbit of
$\rho+\lambda+\sigma_{\alpha\beta}$ is the weight itself. Thus the
only solution to
$[\sigma_{\alpha\beta}+\lambda+\rho]=\lambda_{(ij)}+\rho$ is
$[\sigma_{\alpha\beta}+\lambda+\rho]=[\sigma_{ij}+\lambda+\rho]=\sigma_{ij}+\lambda+\rho$.
Thus the formula for the multiplicity $n_{\lambda_{(ij)}}$ reduces
to $n_{\lambda_{(ij)}}=m'_{\sigma_{ij}}det(Id)=1.1=1.$ $\Box$

We have proved that there is a unique $\G_0$-submodule $\mathrm N$ in
$\lag_2\otimes\V_\lambda\otimes\Sp_\pm$ such that
$\lag_2\otimes\V_\lambda\otimes\Sp_\pm\cong\V_{\lambda_{(ij)}}\otimes\Sp_\pm\oplus\mathrm N$.
Moreover the module $\mathrm N$ is a $\LGP$-submodule of $\mM^\pm_\lambda$.
Let us denote the $\LGP$-module $\mM^\pm_\lambda/\mathrm N$ by
$\U^\pm_\lambda$. 

From now on we will consider the case $\Sp_\pm=\Sp_+$.\footnote{The other case
$\Sp_\pm=\Sp_-$ goes through by simultaneously substituting everywhere $\Sp_-$ for $\Sp_+$ and
$\Sp_+$ for $\Sp_-$.} Let us denote by $\T$ the twistor representation of
$\Spin(n)$. The final piece of information is that
\begin{equation}\label{decomposition of middle piece}
(\V\otimes\E)\otimes(\V_\lambda\otimes\Sp_+)\cong(\V_{\lambda_{(i)}}\otimes\Sp_-)\oplus(\V_{\lambda_{(j)}}\otimes\Sp_-)\oplus(\V_{\lambda_{(i)}}\otimes\T)\oplus(\V_{\lambda_{(j)}}\otimes\T)
\end{equation}
as $\G_0$-modules.

\subsection{Associated Homogeneous Bundle.}
We denote associated vector bundles by the curly letters. A choice of a Weyl structure gives a reduction of the structure group
$\LGP$ of the associated vector bundle $\mathcal
U^+_\lambda=\G\times_\LGP\U^+_\lambda$ to the reductive subgroup
$\G_0$. From the formulas (\ref{g0 pieces in lap}), (\ref{decomposition of big module}), (\ref{decomposition of middle piece}) follows that the bundle $\mathcal U^\pm_\lambda$ decomposes into $\G_0$-subbundles
\begin{eqnarray}
\label{grad on bdle}
\left(
\begin{array}{c}
\mathcal E_0\\
\mathcal V\otimes\mathcal E\\
\Lambda^2\mathcal V
\end{array}
\right)\otimes(\mathcal V_{\lambda}\otimes\mathcal S_+)
\cong\left(
\begin{array}{c}
\mathcal{V}_{\lambda}\otimes\mathcal{S}_{+}\\
(\mathcal{V}_{\lambda_{(i)}}\otimes\mathcal{S}_{-})\oplus
(\mathcal{V}_{\lambda_{(j)}}\otimes\mathcal{S}_{-})\oplus
(\mathcal{V}_{\lambda_{(i)}}\otimes\mathcal{T})\oplus
(\mathcal{V}_{\lambda_{(j)}}\otimes\mathcal{T})\\
\mathcal{V}_{\lambda_{(ij)}}\otimes \mathcal{S}_{+}
\end{array}
\right).
\end{eqnarray}
There is a natural bundle map, which is in the abstract index notation given by
\begin{eqnarray}\label{algebraic action of sections}
&\Gamma
\left(
\begin{array}{c}
\mathcal V\otimes\mathcal E\\
\wedge^2\mathcal V\\
\end{array}
\right)
\otimes\Gamma(\mathcal U^+_\lambda)\ra\Gamma(\mathcal U_\lambda^+)&\\
&\left(
\begin{array}{c}
u_i\otimes\varepsilon_\alpha\\
v_j\wedge w_k 
\end{array}
\right)\otimes\left(
\begin{array}{c}
\varphi_0 \\
x_r\otimes\zeta_\mu\otimes\varphi_1\\
y_s\wedge z_t\otimes\varphi_2
\end{array}
\right)\mapsto\left(
\begin{array}{c}
0\\
u_i\otimes\varepsilon_\alpha\otimes\varphi_0\\
\pi_{ij}(\varepsilon_\alpha\zeta^\alpha u_{[i}x_{r]}\otimes\varphi_1+v_j\wedge w_k\otimes\varphi_0)
\end{array}
\right),
&\nonumber
\end{eqnarray}
where we write section of $\mathcal U_\lambda^+$ as on the left side of (\ref{grad on bdle}), i.e. $u_i,v_j,w_k,x_r,y_s,z_t\in\Gamma(\mathcal V),\varepsilon_\alpha,\zeta_\beta\in\Gamma(\mathcal E),\varphi_m\in\Gamma(\mathcal V_\lambda\otimes\Sp_+)$. The contraction is with respect to the natural scalar product on $\Gamma(\mathcal E)$. The map
\begin{equation}
\pi_{ij}:\Gamma(\Lambda^2\mathcal V\otimes\Sp_+)\ra\Gamma(\mathcal V_{\lambda_{(ij)}}\otimes\Sp_+)
\end{equation}
is the canonical algebraic projection. The pairing (\ref{algebraic action of sections}) comes from the natural action of $\lap_+$ on any finite dimensional $\LGP$-representation. The space of sections $\Gamma(\mathcal U_\lambda^+)$ is thus a module over the Lie algebra $\Gamma(\G\times_\LGP\lap_+)$ and we denote the module structure by $\bullet$.

\section{The Curved Casimir operator.}
We will recall the invariant definition of the Curved Casimir
operator from \cite{CS}.
\begin{df}\label{definition of curved casimir operator}
Let $(\PG\ra M,\omega)$ be a parabolic geometry of type $(\G,\LGP)$.
The bundle $\ATB\M=\PG\times_\LGP\lag$ is called the adjoint tractor
bundle and let $\ATB^\ast\M$ be the dual bundle. Let $\mathcal{WM}$
be the vector bundle induced by a $\LGP$-module $\W$. The Curved
Casimir operator $\mathcal C$ is invariantly defined as the
composition
\begin{equation}\label{invariant def of cco}
\xymatrix{\mathcal C:\Gamma(\mathcal{WM})\ar[r]^{\!\!\!\!\!\!\!\!\!\DO^2}&\Gamma(\otimes^2\ATB^\ast\M\otimes\mathcal W\M)\ar[r]^{\ \ \ \ \ \ B}&\Gamma(\mathcal W\M)},
\end{equation}
where $B$ is the pairing induced on $\ATB\M$ and thus also on $\ATB^\ast\M$ by the Killing form of
$\lag$ and $\DO$ is the fundamental derivative.
\end{df}
The $\lag$-valued Cartan form $\omega$ trivializes the bundle $T\PG$
and identifies sections of $\ATB\M$ with $\LGP$-invariant vector
fields on $\PG$. Let $s\in\Gamma(\ATB\M)$ be a section and let $X$ be the
corresponding vector field on $\PG$. Let $\psi\in\Gamma(WM)$ be a
section of the induced vector bundle as in the definition
\ref{definition of curved casimir operator}. and let $f$ be the
corresponding equivariant function. Then the fundamental derivative
is a pairing
\begin{eqnarray}\label{fundamental derivative}
&\DO:\Gamma(\ATB\M)\otimes\Gamma(\mathcal{WM})\ra\Gamma(\mathcal{WM})&\\
&(s,\psi)\mapsto\DO_s\psi,&\nonumber
\end{eqnarray}
where $\DO_s\psi\in\Gamma(\mathcal{WM})$ is the section which
correspond to the $\LGP$-equivariant function $X.f$ where $X.f$ is
just the derivation of the function $f$ with respect to the vector
field $X$.

Now we give the definition (\ref{invariant def of cco}) of $\mathcal C$ in a local trivialization. Let us denote by $\ATB^i\M=\PG\times_\LGP\lag^i$ the vector subbundles of the adjoint tractor bundle corresponding to the $\LGP$-submodules
$\lag^i=\oplus_{j\ge i}\lag_j$ of the module $\lag$. Let $\mathcal
U\subset M$ be an open set over which the adjoint tractor bundle
$\ATB\M$ is a trivial vector bundle. Let $\ATB\ \mathcal U$ be the
pullback of the adjoint tractor bundle under the inclusion $\mathcal
U\ra M$. Over the set $\mathcal U$, we can choose sections
$s_i,t_j\in\Gamma(\ATB\ \mathcal U),i,j=1,\ldots,\dim(\lag_-)$ such
that:
\begin{enumerate}
\item the sections $s_i$ trivialize the vector bundle $\ATB\ \mathcal U/\ATB^0\mathcal U$.
\item the sections $t_j$ trivialize the vector bundle $\ATB^1\mathcal U$.
\item the sections $s_i,t_j$ are dual with respect to the pairing
$B$, i.e. $B(s_i,t_j)=\delta_{ij}$ at any point $x\in\mathcal U$.
\end{enumerate}
\begin{thm}\label{local description of cco}
Let $\rho\in\Gamma(\mathcal WM)$ and let us use the notation as above. Then on the set $\mathcal U$ we have
\begin{equation}\label{action of cco}
\mathcal C(\rho)|_\mathcal U=-2\sum_{i}t_{i}\bullet\DO_{s_{i}}\rho+c^{\mathfrak{g}_{0}}\rho,
\end{equation}
where $c^{\lag_0}$ is a zero order operator computable from
representation data of the module $\W$ and $\bullet$ denotes the algebraic action of vertical vector fields.
\end{thm}

\subsection{Algebraic action of the Curved Casimir operator.}
The algebraic action $c^{\lag_0}$ of $\mathcal C$ on $\mathcal U_\lambda^+$ is 
\begin{eqnarray}\label{alg action of cco on graded bdle}
&c^{\lag_0}\left(
\begin{array}{ccc}
&\varphi &  \\
&\varphi^{\Sp_-}_i\oplus
\varphi^{\Sp_-}_j\oplus
\varphi^\T_i\oplus
\varphi^\T_j& \\
&\varphi_{ij}&
\end{array}
\right)=\left(
\begin{array}{ccc}
&c_\lambda\varphi &  \\
&c^{\Sp_-}_i\varphi^{\Sp_-}_i\oplus
c^{\Sp_-}_j\varphi^{\Sp_-}_j\oplus
c^\T_i\varphi^\T_i\oplus
c^\T_j\varphi^\T_j& \\
&c_{\lambda_{ij}}\varphi_{ij}&
\end{array}
\right),&
\end{eqnarray}
where we write section with respect to the decomposition on the right side of (\ref{grad on bdle}). In particular $c^{\lag_0}$ acts on each irreducible $\G_0$-subbundle as a multiple of the identity. Formulas for the constants $c_\ast^\bullet$ are explicitly given in the proof of the theorem \ref{gen conf wght}. Let us denote by $\alpha_{\lambda\ast}^\bullet:=c_\lambda-c_{\lambda_\ast}^\bullet$ and $\alpha_{ij}=c_\lambda-c_{\lambda_{ij}}$.

\subsection{Splitting operator from the Curved Casimir operator.}
Let $\pi:\Gamma(\mathcal U_\lambda^+)\ra\Gamma(\mathcal
V_\lambda\otimes\Sp_+)$ be the canonical projection. A splitting
operator $S$ for $\pi$ is a differential operator
\begin{equation}\label{split op}
S:\Gamma(\mathcal V_\lambda\otimes\Sp_+)\ra\Gamma(\mathcal U_\lambda^+)
\end{equation}
such that $\pi\circ S$ is a scalar multiple of the identity operator
on $\Gamma(\mathcal V_\lambda\otimes\Sp_+)$.
\begin{thm}
Let $i:\Gamma(\mathcal V_\lambda\otimes\Sp_+)\ra\Gamma(\mathcal U_\lambda^+)$ be the inclusion of $\G_0$-subbundle given by some Weyl structure Let us denote by $E:=\prod_{\ast,\bullet}(\mathcal C-c_\ast^\bullet)$, where the product is running over all $\G_0$-submodules in $\U_\lambda^+$ except the submodule $\E_0\otimes\V_\lambda\otimes\Sp_+$. The operator $E\circ i$ is a splitting operator for $\pi$.
\end{thm}
Proof: See \cite{CS}.$\Box$

The fact that $S$ is a splitting operator is equivalnt to the following: if we substitute in (\ref{alg action of cco on graded bdle}) $c^{\lag_0}$ by $E$, then the right hand side will depend only on $\varphi$ but not on $\varphi_\ast^\bullet,\varphi_{ij}$. We will use the operator $S:=E\circ i$ to get the operator (\ref{diff op}). It remains to compute the differential part of the Curved Casimir operator.

\section{Affine subset of the homogeneous space $\G/\LGP$.}
Details and proofs, which has been omitted in this paragraph, can be
found in \cite{Ko}. Let $\pi:\G\ra\G/\LGP$ be the homogeneous space
of the parabolic geometry. There is a distribution $\dis\subset T\G/\LGP$ which is the projection of the distribution spanned by the vector fields on $\G$ corresponding to $\lag^1$ under the left trivialization of $T\G$.  This distribution is non-integrable and the Lie bracket is given by the Lie bracket on the Lie algebra $\lag$. Let $\G_-$ be the closed analytic
subgroup $\G_-:=\exp(\lag_-)$ of $\G$. Then $\mathcal U:=\pi(\G_-)$
is an open dense subset of the homogeneous space and $\pi$ is a
diffeomorphism between $\G_-$ and $\mathcal U$. We will denote the $\LGP$-principal fibre bundle $\pi^{-1}(\mathcal U)$ by $\PG$.

\subsection{Preferred trivialization of principle bundles over $\mathcal U$ and flat Weyl structure.}
Let $\mu:\mathcal
U\ra\PG$ be the section determined by
\begin{equation}\label{section of P-bundle}
\mu(\pi(\exp(X)))=\exp(X)
\end{equation}
for all $X\in\lag_-$. Let $\LGP_+:=\exp(\lap_+)$ and let us denote by $\PG_0:=\PG/\LGP_+$. The quotient $\PG_0$ is a principal $\G_0$-bundle over $\mathcal U$. Let $\pi':\PG\ra\PG_0$ be the
canonical projection. Let us denote by $\rho$ the gauge
\begin{eqnarray}\label{trivialization2}
&\varrho:\mathcal U\ra\PG_0&\\
&\varrho(x)=\pi'\circ\mu(x).&\nonumber
\end{eqnarray}

Let $\sigma$ be the Weyl structure, i.e. a $\G_0$-equivariant section, determined by the commutativity of the following diagram
\begin{equation}\label{flat weyl structure}
\xymatrix{&\PG\\
\mathcal U\ar[ru]^\mu\ar[r]^\rho&\PG_0\ar[u]_\sigma.\\}
\end{equation}
The Weyl structure $\sigma$ pulls back the Maurer-Cartan form
$\omega$ to $\PG_0$ and the $\lag_0$-part of $\sigma^\ast\omega$ is principal connection on the $\PG_0$ which we denote by $\omega_0$. The principal connection $\omega_0$ induces
connection on any associated vector bundle and in particular the
induced linear connection on $T\mathcal U$ is flat. Moreover, in the
gauge $\rho$, the connection form, i.e. Christoffel symbols,
$\rho^\ast\omega_0$ is zero.

\subsection{Preferred local adapted frame over $\mathcal U$.}
Let $(\W,\rho)$ be a $\LGP$-module and let $w\in \W$ be a vector.
Let us define a $\LGP$-equivariant function $f_w$ on
$\PG$ by the formula
\begin{eqnarray}\label{equivariant functions}
&f_w:\PG\ra\W&\\
&f_w(g)=\rho^{-1}(p)w,&\nonumber
\end{eqnarray}
where $g\in\PG$ is written uniquely as
$g=\exp(X)p$ with $X\in\lag_-,p\in\LGP$.

In the special case of the adjoint representation $\W=\lag$, we denote by $\hat X$ the $\LGP$-equivariant vector field
on $\PG$ corresponding to the section $f_X\in\Gamma(\ATB\ \mathcal U), X\in\lag$. The value of the vector field $\hat X$ at a point $\exp(Y)p\in\PG$, with
$Y\in\lag_-,p\in\LGP$, is the vector tangent to the curve
$\exp(Y)\exp(tX)p,t\in(-\eps,\eps)$ for $\eps>0$ sufficiently small at $t=0$. Later on we
will need following observations.
\begin{lemma}
Let $X\in\lag_-,Z\in\lag$. Then the derivation of the function $f_Z$
with respect to the vector field $\hat X$ is equal to
\begin{equation}\label{der of eq fun with hor field}
\hat X.f_Z=0
\end{equation}
on $\mathcal U$.
\end{lemma}
Proof:
\begin{eqnarray*}
\hat X.f_Z(\exp(Y)p)=\frac{d}{dt}\big|_0f_Z(\exp(Y)\exp(tX)p)=\frac{d}{dt}\big|_0f_Z(\exp(Y+tX+\frac{1}{2}[Y,tX])p)=0,
\end{eqnarray*}
where we use the Baker-Campbell-Hausdorff formula for $\lag_-$.
$\Box$
\begin{lemma}
Let $Z\in\lag,Z'\in\lap_+$. Then the derivation of the function
$f_Z$ with respect to the vector field $\hat Z'$ is equal to
\begin{equation}\label{der of eq fun with ver field}
\hat Z'.f_Z=f_{-[Z',Z]}=f_{[Z,Z']}.
\end{equation}
\end{lemma}
Proof: As in the formula (\ref{der of eq fun with hor field}) we
have
\begin{eqnarray*}
\hat Z'.f_Z(\exp(Y)p)=\frac{d}{dt}\bigg|_0f_Z(\exp(Y)\exp(tZ')p)=\frac{d}{dt}\bigg|_0p^{-1}e^{-tZ'}Ze^{tZ'}p=f_{-[Z',Z]}=f_{[Z,Z']}.\Box\\
\end{eqnarray*}
Let $p:\ATB\ \mathcal U\ra T\mathcal U$ be the canonical projection.
Let $\xi_X\in\mathfrak{X}(\mathcal U)$ be the vector field
$\xi_X=p_\ast(\hat X)$.
\begin{lemma}\label{bracket of vector fields}
Let $X,Y\in\lag_-$ and let $\xi_X,\xi_Y\in\mathfrak X(\mathcal U)$
be vector fields on the open set $\mathcal U$ defined above. Then
the composition
\begin{eqnarray}\label{left inv vector fields}
&\lag_-\longrightarrow\Gamma(\ATB\ \mathcal U)\longrightarrow\mathfrak X(\mathcal U)&\\
&X\mapsto \hat X\mapsto p(\hat X)&\nonumber
\end{eqnarray}
is a homomorphism of Lie algebras, i.e.
$[\xi_X,\xi_Y]_{\mathfrak{X}(\mathcal
U)}=p\big(\widehat{[X,Y]}_\lag\big)$.
\end{lemma}
Proof: See \cite{Ko}.

\subsection{Functions on $\mathcal U$ and sections.} Let $\W$ be a $\LGP$-representation. Let us use $\mu$ to define an isomorphism
\begin{eqnarray}\label{isom between sections and functions}
&\beta:\mathcal C^\infty(\PG,\W)^\LGP\longrightarrow\mathcal C^\infty(\mathcal U,\W)&\\
&f\longrightarrow f\circ\mu.&\nonumber
\end{eqnarray}
We will write for simplicity $\beta(f)=\tilde f$. The function $f_w$ in (\ref{equivariant functions}) is then the constant function $\tilde f_w(x)=w$ for all $x\in\mathcal U$ and we will denote it for the sake of brevity only as $w$. The inverse map $\beta^{-1}$ we will need only in the simple form $\beta^{-1}(\tilde f)=f$ and in the special case $\tilde f=w$ is constant function then $\beta^{-1}(\tilde f)=f=f_w$.

In particular for any $X\in\lag_-$ and $f\in\mathcal C^\infty(\PG,\W)^\LGP$, we have that
\begin{equation}\label{der of functions and sections}
\beta(\hat X.f)=\xi_X.\tilde f.
\end{equation}
For more see \cite{Ko}.

\subsection{Restriction on sections over $\mathcal U$.}
The freedom of $\LGP$-equivariant functions on $\PG$ is the set $\G_-$. To compare the operators living in the Euclidean and
parabolic setting, we consider functions on the affine subset of the
homogeneous functions which satisfy some additional restrictions. We
will do that in the last section where are given formulas for the
case $k=3$. The case $k=2$ is given without this restriction. Here
are some preliminary notations which will be used later.

Let $\G_2$ be the closed subgroup $\G_{-2}:=\exp(\lag_{-2})$ of
$\G_-$. Let $\tilde{\PG}$ be the quotient of
$\PG$ by the natural left action of $\G_{-2}$. Then
\begin{equation}\label{principal g(-2) bundle}
q:\PG\ra\tilde\PG
\end{equation}
is a principal $\G_{-2}$-bundle. Let us denote by $\tilde{\mathcal
U}$ the quotient space $\G_{-2}\backslash\mathcal U$ by the induced
action of $\G_{-2}$ on $\mathcal U$. The projection $\pi$ is
diffeomorphism between the right coset space $\G_{-2}\backslash\G_-$ and
$\tilde{\mathcal U}$. Moreover $\pi$ descends to $\LGP$-principal
bundle $\tilde{\pi}:\tilde{\PG}\ra\tilde{\mathcal
U}$.

Let $\G_{-1}$ be the subset $\G_{-1}:=\exp(\lag_{-1})$ of $\G_-$.
Let $X\in\lag_-$ be a vector, then there are unique vectors
$X^1\in\lag_{-1},X^2\in\lag_{-2}$  such that $X=X^1+X^2$.
Mapping $\exp(X)\mapsto\exp(X^1)$ gives isomorphism of the right
coset space $\G_{-2}\backslash \G_-$ with $\G_{-1}$ and thus also
identifies $\G_{-1}$ with $\tilde{\mathcal U}$.

\section{Construction of the operator (\ref{diff op}).}
In this section we derive explicit formula for the first term in the theorem \ref{local description of cco}. The operator $S$ given in (\ref{splitting operator}) is a polynomial combination of Curved Casimir operators, in particular it is an operator of degree five in the Curved Casimir operators. However, since the algebraic action of vertical vector fields is compatible with the gradation on the bundle (\ref{grad on bdle}), i.e. the algebraic action of vector fields corresponding to $\lap_+$ raises the homogeneity of sections as can be seen in (\ref{algebraic action of sections}), in the final formula (\ref{final diff operator}) will appear at most second order operator given by differentiating with vector fields in the distribution given by $\lag^1/\lap$. We will need in this section the formulas (\ref{der of eq fun with hor field}),
(\ref{der of eq fun with ver field}) and the fact that $\lap_+$ acts
trivially on irreducible $\LGP$-representations.

Let us first choose a basis $\{X_i|i=1,\ldots,dim(\lag_-)\}$ of $\lag_-$
consisting of homogeneous elements. Let
$\{Z_i|i=1,\ldots,dim(\lap_+)\}$ be the basis of $\lap_+$ dual to
the basis $\{X_i|i=1,\ldots,dim(\lag_-)\}$ with respect to the Killing form of $\lag$. The homogeneity of
elements will be encoded by upper index, for example
$X^1\in\lag_{-1},Z^1\in\lag_1$ etc.

Any $\LGP$-equivariant function of $\PG\times_\LGP\U_\lambda^+$ can be written as $\sum_if_{Z_i}\otimes f_i+f_{e_0}\otimes f_0$. Then
\begin{eqnarray}\label{first use of cco}
&&2\sum_{k=1}^{dim(\lag_-)}\hat Z_k.\hat X_k.(\sum_if_{Z_i}\otimes f_i+f_{e_0}\otimes f_0)\\
&=&2\sum_{i,k}(\hat Z_k.f_{Z_i}\otimes\hat X_k.f_i+\hat Z_k.f_{e_0}\otimes\hat X_k.f_0)\nonumber\\
&=&-2\sum_{i,k}(f_{[Z_k,Z_i]}\otimes\hat X_k.f_i+f_{[Z_k,e_0]}\otimes\hat X_k.f_0)\nonumber
\end{eqnarray}

In the gauge $\mu$, the formula (\ref{first use of cco}) is
\begin{equation}\label{first use of cco in triv}
2\sum_{k=1}^{dim(\lag_-)}Z_k.\xi_{X_k}.
\left(
\begin{array}{c}
e_0\otimes\tilde f_0\\
\sum_i Z_i^1\otimes\tilde f_i^1\\
\sum_j Z_j^2\otimes\tilde f_j^2
\end{array}
\right)=\left(
\begin{array}{c}
0\\
2\sum_kZ^1_k\otimes\xi_{X^1_k}\tilde f_0 \\
\pi_{ij}(4\sum_kZ^2_k\otimes\xi_{X^2_k}\tilde f_0-2\sum_{ki}[Z_k^1,Z^1_i]\otimes\xi_{X^1_k}.\tilde f_i^1)\\
\end{array}
\right).
\end{equation}

Let us set $D^1:=\sum_{i=1}^{dim(\lag_1)}Z^1_i.\xi_{X^1_i}$. Then $D^1$ is first order differential operator on the set $\mathcal U$.
Let us denote by $D^\W_\ast$ the second order operator on $\mathcal U$ given by
\begin{eqnarray}\label{second order operators}
\xymatrix{D^\bullet_\ast:\mathcal C^\infty(\mathcal U,\V_\lambda\otimes\Sp_+)\ar[rr]^{\pi^\W_\ast\circ D^1}&&\mathcal C^\infty(\mathcal U,\V_{\lambda_\ast}\otimes\W)\ar[rr]^{\pi\circ D^1}&&\mathcal C^\infty(\mathcal U,\V_{\lambda_{(ij)}}\otimes\Sp_+)}
\end{eqnarray}
where $\ast\in\{i,j\}$ and $\W$ stands for the representation $\Sp_-$ or $\T$. The symbol $\pi^\W_\ast$ is the algebraic projection $\pi^\W_\ast:(\V\otimes\E)\otimes(\V_\lambda\otimes\Sp_-)\ra\V_{\lambda_\ast}\otimes\W$. In particular, the second map in (\ref{second order operators}) is the composition
\begin{eqnarray}\label{second par second order op}
\mathcal C^\infty(\mathcal U,\V_{\lambda_\ast}\otimes\W)\longrightarrow^{\!\!\!\!\!\!\!\!\!D^1}\mathcal C^\infty(\mathcal U,(\V\otimes\E)\otimes(\V\otimes\E)\otimes\V_{\lambda}\otimes\Sp_+)\longrightarrow^{\!\!\!\!\!\!\!\!\!\rho}\\
\longrightarrow^{\!\!\!\!\!\!\!\rho}\ \mathcal C^\infty(\mathcal U,(\Lambda^2\V\otimes\C)\otimes(\V_\lambda\otimes\Sp_+))\longrightarrow^{\!\!\!\!\!\!\!\!\!\!\!\pi_{ij}}\ \mathcal C^\infty(\mathcal U,\V_{\lambda_{(ij)}}\otimes\Sp_+).\nonumber
\end{eqnarray}
In particular $\rho=(\wedge\otimes Tr)\otimes Id_{\V_\lambda\otimes\Sp_+}$ is a tensor product of the Lie bracket $\wedge\otimes Tr$ given in (\ref{lie bracket}) with the identity map on $\V_\lambda\otimes\Sp_+$. The map $\pi_{ij}$, as in the formula (\ref{algebraic action of sections}), is the natural projection.

With this notation we can give formula for the operator
\begin{equation}\label{splitting operator}
E:=(\mathcal{C}-c_{\lambda_{i}^\Sp})(\mathcal{C}-c_{\lambda_{j}^\Sp})
(\mathcal{C}-c_{\lambda_{i}^\T})(\mathcal{C}-c_{\lambda_{j}^\T})(\mathcal{C}-c_{\lambda_{ij}})
\end{equation}
introduced below the formula (\ref{split op}). We have that
\begin{eqnarray}\label{final diff operator}
&&(\mathcal{C}-c_{\lambda_{i}^\Sp})(\mathcal{C}-c_{\lambda_{j}^\Sp})
(\mathcal{C}-c_{\lambda_{i}^\T})(\mathcal{C}-c_{\lambda_{j}^\T})(\mathcal{C}-c_{\lambda_{ij}})
\left(\begin{array}{c}
e_0\otimes\tilde f\\
\ast\oplus\ast\oplus\ast\oplus\ast\\
\ast\\
\end{array}
\right)\nonumber\\
&=&(\mathcal{C}-c_{\lambda_{j}^\Sp})(\mathcal{C}-c_{\lambda_{i}^\T})(\mathcal{C}-c_{\lambda_{j}^\T})(\mathcal{C}-c_{\lambda_{ij}})
\left(\begin{array}{c}
\alpha_i^\Sp\otimes\tilde f\\
2\pi_i^\Sp(D^1f)\oplus\ast\oplus\ast\oplus\ast\\
\ast\\
\end{array}
\right)\nonumber\\
&=&(\mathcal{C}-c_{\lambda_{i}^\T})(\mathcal{C}-c_{\lambda_{j}^\T})(\mathcal{C}-c_{\lambda_{ij}})
\left(\begin{array}{c}
\alpha_i^\Sp\alpha_j^\Sp e_0\otimes\tilde f\\
2\alpha_j^\Sp\pi_i^\Sp(D^1 f)\oplus2\alpha_i^\Sp\pi_i^\Sp(D^1 f)\oplus\ast\oplus\ast\\
\ast\\
\end{array}
\right)\nonumber\\
&=&\ldots\nonumber\\
&=&(\mathcal{C}-c_{\lambda_{ij}})
\left(\begin{array}{c}
\alpha_i^\Sp\alpha_j^\Sp\alpha_i^\T\alpha_j^\T e_0\otimes\tilde f\\
2(c_i^\Sp\pi_i^\Sp\oplus c_j^\Sp\pi_j^\Sp\oplus c_i^\T\pi_i^\T\oplus c_j^\T\pi_j^\T)(D^1 f)\\
\ast\\
\end{array}
\right)
\nonumber\\
&=&\left(
\begin{array}{c}
\ \ \ ae_0\otimes\tilde f\\
2(b_i^\Sp\pi_i^\Sp\oplus b_j^\Sp\pi_j^\Sp\oplus b^\T_i\pi^\T_i\oplus b^\T_j\pi^\T_j)(Z^1_l\otimes\xi_{X^1_l}.\tilde f)\\
4\pi_{ij}\big(c_i^\Sp D^\Sp_i\tilde f+c_j^\Sp D^\Sp_j\tilde f+c_i^\T D^\T_i\tilde f+c_j^\T D^\T_j \tilde f+c_{ij}\sum_lZ^2_l\otimes\xi_{X^2_l}.\tilde f\big)\\
\end{array}
\right),
\end{eqnarray}
where the coefficients are
\begin{eqnarray*}\label{coefficients}
a=\alpha_j^\Sp\alpha_j^\Sp\alpha_i^\T\alpha_i^\T\alpha_{ij},
b_\ast^\bullet=\frac{a}{\alpha_\ast^\bullet},
c_\ast^\bullet=\frac{a}{\alpha_{ij}\alpha_\ast^\bullet},c_{ij}=\frac{a}{\alpha_{ij}}.
\end{eqnarray*}
The coefficients $\alpha_\ast^\bullet$ are given below (\ref{alg action of cco on graded bdle}). Thus, if
\begin{equation}\label{con eq}
\alpha_{ij}=0,
\end{equation}
then the first two rows in the formula (\ref{final diff operator}) are zero. The Curved Casimir operator is $\G$-invariant and thus also the operator in (\ref{final diff operator}) is also $\G$-invariant. Let us summarize it into the following theorem.
\begin{thm}\label{1. thm}
Suppose that (\ref{con eq}) holds. Then the second order differential operator $D$ given in the last row in the formula (\ref{final diff operator}) is $\G$-invariant operator
\begin{eqnarray}\label{diff op from splitting op}
D:\Gamma(\mathcal V_\lambda\otimes\mathcal{S}_+)\longrightarrow\Gamma(\mathcal V_{\lambda_{(ij)}}\otimes\mathcal S_+),
\end{eqnarray}
as in (\ref{diff op}).
\end{thm}
The equation (\ref{con eq})\ is an equation for generalized
conformal weight. Now we show that (\ref{con eq}) holds. In the last
two sections I give the formulas of the operators coming from
(\ref{diff op from splitting op}) with the coefficients in the
theorem \ref{wght 1}.
\begin{thm}\label{gen conf wght}
Let $n\ge 2k\ge 4$ and suppose that $n$ is even. Let $\V_\lambda\otimes\Sp_\pm,\V_\nu\otimes\Sp_\pm$ be two
irreducible $\LGP$-modules such that the highest weights of the dual
modules lie either on the affine orbit of the weight
$\frac{1}{2}(1-n,\ldots,1-n|1,\ldots,1)$ or on the affine orbit of
the weight $\frac{1}{2}(1-n,\ldots,1-n|1,\ldots,1,-1)$ and moreover
suppose that $\nu=\lambda_{(ij)}$ for some $1\le i<j\le k$. Then the equation (\ref{con eq}) holds.
\end{thm}
Proof: Let $\W$ be an irreducible $\G_0$-module with the lowest
weight $-\mu$, the highest weight $\upsilon$ and let $\delta$ be the
lowest weight of $\mathfrak{g}$. The (algebraic) action of the
Curved Casimir operator on the sections of the bundle induced by
$\W$ is equal to $\langle\mu,\mu+2\delta\rangle$. Then for the
module $\W$ we have that
\begin{equation*}
\upsilon=\big(\mu_1,\ldots,\mu_k\big|\mu_{k+1},\ldots,\mu_{k+n}\big)\leftrightarrow
-\mu=\big(\mu_k,\ldots,\mu_1|-\mu_{k+1},\ldots,-\mu_{k+n}\big).
\end{equation*}
Thus the difference $c_{\lambda}-c_{\lambda_{ij}}$ is equal to
\begin{eqnarray}\label{help2}
(n+2k-2)(c_{\lambda}-c_{\lambda_{ij}})&=&-\lambda_{i}(-\lambda_{i}+2(\frac{n}{2}+i-1))-\lambda_{j}(-\lambda_{j}+2(\frac{n}{2}+j-1))\nonumber \\
&-&\big[(-\lambda_{i}-1)(-\lambda_{i}-1+2(\frac{n}{2}+i-1))-(\lambda_{j}-1)(-\lambda_{j}-1+2(\frac{n}{2}+j-1))\big]\nonumber\\
&=&-2\lambda_{i}-1+2(\frac{n}{2}+i-1)-2\lambda_{j}-1+2(\frac{n}{2}+j-1))\nonumber\\
&=&-2\lambda_{i}-2\lambda_{j}+2n+2i+2j-6.
\end{eqnarray}
The weight $\lambda$ can be written as
$\lambda=\lambda'+(\lambda-\lambda')$ where
$\lambda'=\frac{1}{2}(n-1,\ldots,n-1)$ and moreover the partition
$\lambda-\lambda'$ has the Young diagram symmetric with respect
to the main diagonal. From the symmetry we get that
\begin{equation}\label{help1}
\lambda_j-\frac{1}{2}(n-1)=i-1;\lambda_i-\frac{1}{2}(n-1)=j-1.
\end{equation}
Plugging (\ref{help1}) into (\ref{help2}) we get that
\begin{eqnarray*}
\alpha_{i,j}=c_{\lambda}-c_{\lambda_{ij}}=\frac{2n-2(n-1)-2}{n+2k-2}=0.
\end{eqnarray*}
\begin{thm}\label{wght 1} The coefficients of the operator (\ref{final diff operator}) are
equal to
$\alpha_i^\Sp=\frac{2(\lambda_j-\lambda_i)}{n+2k-2},\alpha_j^\Sp=\frac{2(\lambda_i-\lambda_j)}{n+2k-2},\alpha_i^\T=\frac{2(\lambda_j-\lambda_i)-n}{n+2k-2},\alpha_j^\Sp=\frac{2(\lambda_i-\lambda_j)-n}{n+2k-2}$.
\end{thm}
Proof: Let us compute for example $c_{\lambda}-c^\Sp_{\lambda_{i}}$.
Then we have
\begin{eqnarray*}
(n+2k-2)c_{\lambda}-c^\Sp_{\lambda_{i}}&=&-2\lambda_{i}-1+2(\frac{n}{2}+i-1)=-2\lambda_{i}+n+2i-3\nonumber\\
&=&-2\lambda_{i}+(n-1)+2i-2=-2\lambda_{i}+2\lambda_{j}=2(\lambda_{j}-\lambda_{i}),
\end{eqnarray*}
where we have used (\ref{help1}). Similarly we get
$c_{\lambda}-c^\Sp_{\lambda_{j}}=2(\lambda_{i}-\lambda_{j})$ and
similarly for the remaining coefficients.

\section{Local formulas of the operators.} Let $\{e_1,e_2,\ldots,e_k\}$ be
the standard basis of $\R^k$ and let $\{e^1,e^2,\ldots,e^k\}$ be the
dual basis of the $\SL(k,\R)$-module $(\R^k)^\ast$. Let
$\{\varepsilon_\alpha,1\le\alpha\le n\}$ be an orthonormal basis of
$\R^n$. We denote the $\so(n)$-invariant product on $\R^n$ by $g_{\alpha\beta}$. Let $B$ the $\G$-invariant scalar product on $\lag$ as in
the remark 3.

The section $(\ref{trivialization2})$ gives an isomorphism
$\phi:\lag_{-1}\cong(\R^k)^\ast\otimes\R^n$. Let
$\{X_{i\alpha}|\sqrt{n+2}X_{i\alpha}=\phi^{-1}(e^i\otimes\varepsilon_\alpha),i=1,\ldots,k,\alpha=1,\ldots,n\}$
be a preferred basis of $\lag_{-1}$ and let
$Z_{j\beta},j=1,2,\beta=1,\ldots,n$ be the dual elements in $\lag_1$
with respect to $B$. Then
$\{X_{ij}=-X_{ji}\big|\sqrt{n+2}[X_{i\alpha},X_{j\beta}]=\delta_{\alpha\beta}X_{ij},1\le i,j\le k,\alpha,\beta\le n\}$ is a preferred basis of $\lag_{-2}\cong\Lambda^2(\R^\ast)^k$. Let
$\{Z_{ij}|1\le i,j\le k\}$ be the basis of $\lag_2$ dual with
respect to $B$ to the basis $\{X_{ij}|1\le i,j\le k\}$ of $\lag_{-2}$, then we have that
$\sqrt{n+2}[Z_{i\alpha},Z_{j\beta}]=-\delta_{\alpha\beta}Z_{ij}$.

Let us write the canonical coordinates on $\lag_-$ given by the preferred basis $\{X_{i\alpha},X_{ij}\}$ by $(x_{i\alpha},x_{ij})$. We may use these coordinates also on the set $\osu$ and let $\partial_{i\alpha}$ and $\partial_{ij}$ be the coordinate vector fields. The left invariant vector fields are then
\begin{eqnarray}\label{left invariant vector fields in the coordinates}
\xi_{X_{k\mu}}(x_{i\alpha},x_{ln})&=&\partial_{k\mu}-\frac{1}{2\sqrt{n+2}}\sum_ix_{i\mu}\partial_{ki}\\
\xi_{X_{rs}}(x_{i\alpha},x_{ln})&=&\partial_{rs}\nonumber.
\end{eqnarray}
The left invariant vector fields $\xi_{X_{k\nu}}$ span the distribution $\dis$ on $\mathcal U$. With this notation, the commutator is
\begin{eqnarray}\label{commutator of left inv fields in coordinates}
[\xi_{X_{j\mu}},\xi_{X_{k\nu}}](x_{i\alpha},y_{mn})&=&(\partial_{j\mu}-\frac{1}{2\sqrt{n+2}}\sum_ix_{i\mu}\partial_{ji})(\partial_{k\nu}- \frac{1}{2\sqrt{n+2}}\sum\limits_ix_{i\nu}\partial_{ki})\nonumber\\
&-&(\partial_{k\nu}- \frac{1}{2\sqrt{n+2}}\sum\limits_ix_{i\nu}\partial_{ki})(\partial_{j\mu}-\frac{1}{2\sqrt{n+2}}\sum_ix_{i\mu}\partial_{ji})\nonumber\\
&=&g_{\mu\nu}\bigg(\frac{-1}{2\sqrt{n+2}}\partial_{kj}+\frac{1}{2\sqrt{n+2}}\partial_{jk}\bigg)=\frac{g_{\mu\nu}}{\sqrt{n+2}}\partial_{jk}.
\end{eqnarray}

For $i=1,\ldots,k$, let us denote by $\partial_i$ the
first order differential operators
$\partial_i=\sum_{\alpha=1}^n\varepsilon_\alpha\xi_{X_{i\alpha}}$
where $\varepsilon_\alpha$ denotes the multiplication with the
Clifford number $\varepsilon_\alpha$. These operators are analogues of the $k$-Dirac operators defined in (\ref{k-dirac operator}) in the parabolic setting.

\subsection{The sequence for $k=2$.}
A regular parabolic geometry of this type is a contact geometry, in
particular $\lag_{-2}$ is one-dimensional. There are three
differential operators in the sequence starting with the $2$-Dirac
operator. The first operator is the $2$-Dirac operator, which we
denote by $D_1$, the second one is a second order operator
\begin{equation}\label{first operator}
D_2:\Gamma(\mathcal V_{\lambda_2}\otimes\mathcal{S}_+)\ra\Gamma(\mathcal
V_{\lambda_3}\otimes\mathcal S_+)
\end{equation}
where the weights are
$\lambda_2=\frac{1}{2}(n+1,n-1),\lambda_3=\frac{1}{2}(n+3,n+1)$ and
the sequence closes with a first order operator $D_3$. The operators
$D^1,D^2$ can be computed easily using the Curved Casimir operator.
Here we will focus on the operator (\ref{first operator}). 

Before giving the formula for the operator, let us first give more explicitly the formula (\ref{algebraic action of sections}). The projections are
\begin{eqnarray}\label{twistor and spinor projections}
\E\otimes\Sp_+&\ra&\T\oplus\Sp_-\\
\varepsilon\otimes\varphi&\mapsto&(\varepsilon\otimes\varphi+\frac{1}{n}\sum_\alpha\varepsilon_\alpha\otimes\varepsilon_\alpha\varepsilon\varphi)\oplus(-\frac{1}{n}\sum_\alpha\varepsilon_\alpha\otimes\varepsilon_\alpha\varepsilon\varphi),\nonumber
\end{eqnarray}
where the sums are running over the standard basis of the defining module $\E$ and $\varepsilon\in\E,\varphi\in\Sp_+$. The sign is coming from the defining relation of the Clifford algebra, i.e. $\varepsilon_\alpha\varepsilon_\beta+\varepsilon_\beta\varepsilon_\alpha=-2g_{\alpha\beta}$. 
Let us work in the gauge $\mu$ as in (\ref{section of P-bundle}). Let $\varphi_\ast$ be spinor valued function on $\mathcal U$ and let us keep the notation introduced below the formula (\ref{isom between sections and functions}). The formula (\ref{algebraic action of sections}) is in this case 
\begin{eqnarray}
&-\left(
\begin{array}{c}
\sum_{s\beta}e_s\otimes\varepsilon_\beta\\
e_1\wedge e_2
\end{array}
\right)\bullet
\left(
\begin{array}{c}
\sum_i e_i\otimes\varphi_i \\
\sum_{jk\mu}e_j\otimes e_k\otimes\varepsilon_\mu\otimes\varphi_{jk\mu}\\
\sum_l e_1\wedge e_2\otimes e_l\otimes\varphi_l
\end{array}
\right)&\nonumber\\
&=\left(
\begin{array}{c}
0\\
\sum_{is\beta\alpha}e_s\wedge e_i\otimes(\varepsilon_\beta\otimes\varphi+\frac{1}{n}\varepsilon_\alpha\otimes\varepsilon_\alpha\varepsilon_\beta\varphi_i)\oplus\sum_{is\beta\alpha}e_s\wedge e_i\otimes(-\frac{1}{n}\varepsilon_\alpha\otimes\varepsilon_\alpha\varepsilon_\beta\varphi_i)\\
\oplus\sum_{is\beta\alpha}e_s\odot e_i\otimes(\varepsilon_\beta\otimes\varphi+\frac{1}{n}\varepsilon_\alpha\otimes\varepsilon_\alpha\varepsilon_\beta\varphi_i)\oplus\sum_{is\beta\alpha}e_s\odot e_i\otimes(-\frac{1}{n}\varepsilon_\alpha\otimes\varepsilon_\alpha\varepsilon_\beta\varphi_i)\\
2\sum_i e_1\wedge e_2\otimes e_i\otimes\varphi_i-\frac{g_{\beta\mu}}{\sqrt{n+2}}\sum_{sjk\mu\beta}e_s\wedge e_j\otimes e_k\otimes\varphi_{jk\mu}\\
\end{array}
\right).&\nonumber
\end{eqnarray}

With all these preliminary results we can finally give the simplest second order operator explicitly.
The operator (\ref{first operator}) is
\begin{equation}\label{k=2 operator}
\left(
\begin{array}{c}
\e_1\otimes\phi_1 \\
\e_2\otimes\phi_2 \\
\end{array}\right)
\mapsto \left(
\begin{array}{c}
\e_1\wedge\e_2\otimes\e_1\otimes(\partial_1\partial_1\phi_2-\partial_2\partial_1\phi_1+\frac{2}{\sqrt{n+2}}\xi_{X_{12}}\phi_1)\\
\e_1\wedge\e_2\otimes\e_2\otimes(\partial_1\partial_2\phi_2-\partial_2\partial_2\phi_1+\frac{2}{\sqrt{n+2}}\xi_{X_{12}}\phi_2)\\
\end{array}
\right).
\end{equation}
The formula (\ref{k=2 operator}) remains unchanged if we swap
$\Sp_+$ and $\Sp_-$. Since the operator (\ref{first operator}) with
$\Sp_-$ is invariant for the same generalized conformal weight we
can take $\Sp^\pm=\Sp_+\oplus\Sp_-$instead of just $\Sp_+$. We do
this in the remaining paragraphs.

\subsection{The sequence for $k=2$ is a complex.}
The composition $D_2\circ D_1$ is equal to
\begin{equation}\label{complex verification}
\phi\mapsto
\left(
\begin{array}{c}
e_1\otimes\partial_1\phi \\
e_2\otimes\partial_2\phi \\
\end{array}
\right)\mapsto
\left(
\begin{array}{c}
\e_1\wedge\e_2\otimes\e_1\otimes(\partial_1\partial_1\partial_2\phi-\partial_2\partial_1\partial_1\phi+\frac{2}{\sqrt{n+2}}\xi_{X_{12}}\partial_1\phi)\\
\e_1\wedge\e_2\otimes\e_2\otimes(\partial_1\partial_2\partial_2\phi-\partial_2\partial_2\partial_1\phi+\frac{2}{\sqrt{n+2}}\xi_{X_{12}}\partial_2\phi)\\
\end{array}
\right).
\end{equation}
From the lemma $6.$ follows that
\begin{equation}\label{help for complex}
\partial_1\partial_1\partial_2\phi-\partial_2\partial_1\partial_1\phi=-\frac{2}{\sqrt{n+2}}\xi_{X_{12}}\partial_1\phi,
\ \ \ \partial_1\partial_2\partial_2\phi-\partial_2\partial_2\partial_1\phi=-\frac{2}{\sqrt{n+2}}\xi_{X_{12}}\partial_2\phi.\\
\end{equation}
Plugging (\ref{help for complex}) into (\ref{complex verification})
gives
\begin{equation}\label{final complex}
\frac{2}{\sqrt{n+2}}\left(
\begin{array}{c}
\e_1\wedge\e_2\otimes\e_1\otimes(-\xi_{X_{12}}\partial_1\phi+\xi_{X_{12}}\partial_1\phi)\\
\e_1\wedge\e_2\otimes\e_2\otimes(-\xi_{X_{12}}\partial_2\phi+\xi_{X_{12}}\partial_2\phi)\\
\end{array}
\right)=
\left(
\begin{array}{c}
0 \\
0 \\
\end{array}
\right).
\end{equation}
The composition of the last two operators is equal to
\begin{eqnarray*}\label{last composition}
&\left(
\begin{array}{c}
e_1\otimes\phi_1 \\
e_2\otimes\phi_2 \\
\end{array}
\right)
\mapsto
\left(
\begin{array}{c}
\e_1\wedge\e_2\otimes\e_1\otimes(\partial_1\partial_1\phi_2-\partial_2\partial_1\phi_1+\frac{2}{\sqrt{n+2}}\xi_{X_{12}}\phi_1)\\
\e_1\wedge\e_2\otimes\e_2\otimes(\partial_1\partial_2\phi_2-\partial_2\partial_2\phi_1+\frac{2}{\sqrt{n+2}}\xi_{X_{12}}\phi_2) \\
\end{array}
\right)\mapsto&\\
\mapsto&\frac{2}{\sqrt{n+2}}\e_1\wedge\e_2\otimes\e_1\wedge
e_2\otimes\big[\partial_1(\partial_1\partial_2\phi_2-\partial_2\partial_2\phi_1+\frac{2}{\sqrt{n+2}}\xi_{X_{12}}\phi_2)\nonumber&\\
&-\partial_2(\partial_1\partial_1\phi_2-\partial_2\partial_1\phi_1+\frac{2}{\sqrt{n+2}}\xi_{X_{12}}\phi_1)\big]&
\nonumber
\end{eqnarray*}
The equations for $\phi_1$ and for $\phi_2$ in (\ref{last composition}) are identical to (\ref{final complex}).
Thus the sequence for $k=2$ is a complex.

\subsection{The symbol sequence for $k=2$.}
Let us recall that we have denoted by $\dis$ the distribution on $T\mathcal U$ corresponding to $\lag_{-1}$. The highest order parts of the operators $D_1,D_2,D_3$ belong to $\dis$ so the symbol of these operators is determined by its restriction to $\dis$. We have that $\dis^\ast$ is a quotient of $T^\ast\mathcal U$. Let $x\in\mathcal U,v\in\dis^\ast_x$. Then for the vector $v$ the symbol sequence is
\begin{equation}\label{symbol sequence}
\xymatrix{
\V_{\lambda_1}\otimes\Sp^\pm\ar[rr]^{\sigma_{(x,v)}(D_1)}&&\V_{\lambda_2}\otimes\Sp^\pm\ar[rr]^{\sigma_{(x,v)}(D_2)}&&\V_{\lambda_3}\otimes\Sp^\pm\ar[rr]^{\sigma_{(x,v)}(D_3)}&&\V_{\lambda_4}\otimes\Sp^\pm}
\end{equation}
where $\lambda_2,\lambda_3$ are given below (\ref{first operator})
and $\lambda_1=\frac{1}{2}(n-1,n-1),\lambda_4=\frac{1}{2}(n+3,n+3)$.
Since the sequence of operators is complex also the sequence
(\ref{symbol sequence}) is a complex. Thus it suffices to show that,
at each point in the sequence, the dimension of the image has
maximal possible dimension. We can consider the symbol map up to a
scalar multiple.

Let us denote $\sum_{\alpha}(\xi_{
X_{i\alpha}}.f)(x)e_i\otimes\epsilon_\alpha$ for $x\in\mathcal U$ by
$f_i\in\dis_x^\ast\cong\V\otimes\E\subset\V\otimes End(\Sp^\pm)$. The symbol of the first
operator $D^1$ is
\begin{equation}
\sigma_{(x,v)}(D_1)=\left(
\begin{array}{c}
f_1 \\
f_2 \\
\end{array}
\right)
\end{equation}
which is an injective map since if $f_i\ne 0$ then $f_i$ is an
isomorphism of spinor spaces. The symbol $\sigma_{(x,v)}(D^2)$ is
\begin{equation}
\sigma_{(x,v)}(D_2)=\left(\begin{array}{cc}
-f_2f_1& f_1^2 \\
-f_2^2 & f_2f_1  \\
\end{array}
\right)
\end{equation}
If $f_i\ne 0$ then $f_i^2$ is a multiple of the identity matrix
$I\in End(\Sp^\pm)$ and thus the rank of the matrix is equal to $dim(\Sp^\pm)$. The symbol of the last operator $D_3$ is
\begin{equation}
\sigma_{(x,v)}(D_3)=\left(
\begin{array}{cc}
-f_2& f_1
\end{array}
\right)
\end{equation}
and so the last symbol is surjective when restricted to $\dis^\ast$.

Although we have computed symbols and checked that the sequence
is a complex for a particular choice of Weyl structure over the
open set $\mathcal U$, both statements remain true if we choose
different Weyl structure over the set $\mathcal U$. Instead of the
set $\mathcal U$, we can consider open covering
$\{\pi(\exp(\lag_-)g)|g\in\G\}$ of the homogeneous space. All
constructions work over these sets as well and thus both statements
hold on $\G/\LGP$.

\begin{thm}\label{complex for k=2}
The sequence of operators for $k=2$ is a complex and the symbol sequence restricted to $\dis$ is exact.
\end{thm}

\section{Sequence for $k=3$.} In this section we will consider smooth $\LGP$-equivarint
function on $\pi^{-1}(\mathcal{U})$ which are in the image of the
pullback $q^\ast:\mathcal C^\infty(\pi^{-1}(\tilde{\mathcal
U}),\W)^\LGP\ra\mathcal C^\infty(\pi^{-1}(\mathcal U),\W)^\LGP$
where $q$ is as in (\ref{principal g(-2) bundle}). Any
$\LGP$-equivariant function on $\pi^{-1}(\mathcal U)$ is determined
by its values on $\G_-$ and any smooth function $f\in Im(q^\ast)$ is
determined by its values on $\G_{-1}$. A smooth function $f$ is in
the image $Im(q^\ast)$ of $q^\ast$ iff $f$ is constant on the orbits
under the right action of $\G_{-2}$ on $\pi^{-1}(\mathcal U)$. Using
the flat Weyl structure (\ref{flat weyl structure}), we see that
this is equivalent to $\xi_Xf=0$ for all $X\in\lag_{-2}$. Since
$\G_{-1}\cong\R^{n+k}$ we have an analogy with the operators
appearing in the resolutions of the operator (\ref{k-dirac
operator}). The sequence for $k=3$ looks like
\begin{equation}
\xymatrix{\Gamma(\mathcal{V}_{\lambda_1}\otimes\mathcal{S}_\pm)\ar[r]^{D_1}&\Gamma((\mathcal{V}_{\lambda_2}\otimes\mathcal{S}_\pm)\ar[r]^{D_2}&\Gamma((\mathcal{V}_{\lambda_3}\otimes\mathcal{S}_\pm))\ar[d]^{D_4}\ar[r]^{D_3}&\Gamma((\mathcal{V}_{\lambda_4}\otimes\mathcal{S}_\pm)\ar[d]^{D_6}&\Gamma((\mathcal{V}_{\lambda_8}\otimes\mathcal{S}_\pm)\\
&&\Gamma(\mathcal{V}_{\lambda_5}\otimes\mathcal{S}_\pm)\ar[r]^{D_5}&\Gamma((\mathcal{V}_{\lambda_6}\otimes\mathcal{S}_\pm)\ar[r]^{D_7}&\Gamma((\mathcal{V}_{\lambda_7}\otimes\mathcal{S}_\pm))\ar[u]^{D_8}}
\end{equation}
where
\begin{eqnarray*}
\lambda_1=\frac{1}{2}(n-1,n-1,n-1)\\
\lambda_2=\frac{1}{2}(n+1,n-1,n-1)\\
\lambda_3=\frac{1}{2}(n+3,n+1,n-1)\\
\lambda_4=\frac{1}{2}(n+3,n+3,n-1)\\
\lambda_5=\frac{1}{2}(n+5,n+1,n+1)\\
\lambda_6=\frac{1}{2}(n+5,n+3,n+1)\\
\lambda_7=\frac{1}{2}(n+5,n+5,n+3)\\
\lambda_8=\frac{1}{2}(n+5,n+5,n+5).\\
\end{eqnarray*}
The operators $D_2,D_4,D_6,D_7$ are second order operators, others
are first order operators. Considering only function belonging to
the $Im(q^\ast)$, one can verify directly that the
sequence is complex, i.e. in the middle box one can show that
$D_6\circ D_3+D_5\circ D_4=0$, and that the symbol sequence is exact on $\dis$.
\begin{thm}\label{complex for k=3}.
The sequence of operators for $k=2$ is a complex and the symbol sequence restricted to $\dis$ is exact.
\end{thm}

Now we will give explicit formulas of the operators. Let
$\{e_1,e_2,e_3\}$ be a basis of $V_{\lambda_2}$. Then
\begin{equation}
D_1(\phi)=
\left(\begin{matrix}
e_1\otimes\partial_1\phi\\
e_2\otimes\partial_2\phi\\
e_3\otimes\partial_3\phi\\
\end{matrix}\right).
\end{equation}
Let
\begin{eqnarray*}
&h_{ij}=e_i\wedge e_j\otimes e_i,i\ne j,1\le i,j\le 3,&\\
&w_1=\frac{1}{3}e_1\wedge e_2\wedge e_3-e_3\wedge e_1e_2,&\\
&w_2=\frac{1}{3}e_1\wedge e_2\wedge e_3-e_2\wedge e_3e_1&
\end{eqnarray*} be a basis of $\V_{\lambda_3}$. With respect to this basis the operator $D^2$ is equal to
\begin{equation}
\left(\begin{array}{c}
e_1\otimes\phi_1\\
e_2\otimes\phi_2\\
e_3\otimes\phi_3\\
\end{array}\right)
\mapsto
\left(\begin{array}{c}
h_{ij}\otimes(\partial_i\partial_i\phi_j-\partial_j\partial_i\phi_i)\\
w_1\otimes(\{\partial_2,\partial_3\}\phi_1-\partial_1\partial_3\phi_2-\partial_1\partial_2\phi_3)\\
w_2\otimes(\{\partial_1,\partial_3\}\phi_2-\partial_2\partial_3\phi_1-\partial_2\partial_1\phi_3)\\
\end{array}\right).
\end{equation}
For brevity we omit the $\otimes$ symbol. Let
\begin{eqnarray*}
&&v_{220}=e_1e_1\wedge e_2e_2-e_2e_1\wedge e_2e_1,\\
&&v_{202}=e_1e_1\wedge e_3e_3-e_3e_1\wedge e_3e_1,\\
&&v_{022}=e_2e_2\wedge e_3e_3-e_3e_2\wedge e_3e_2,\\
&&v_{211}=e_1e_1\wedge e_3e_2+e_1e_1\wedge e_2e_3-e_3e_1\wedge e_2e_1-e_2e_1\wedge e_3e_1,\\
&&v_{121}=e_3e_1\wedge e_2e_2+e_1e_3\wedge e_2e_2-e_2e_3\wedge e_2e_1-e_2e_1\wedge e_2e_3,\\
&&v_{112}=e_3e_1\wedge e_3e_2+e_1e_3\wedge e_2e_3-e_3e_3\wedge e_2e_1-e_2e_1\wedge e_3e_3\\
\end{eqnarray*} be a basis of $\V_{\lambda_4}$. With respect to this basis the operator $D_3$ is equal to
\begin{equation}
\left(\begin{array}{c}
h_{ij}\otimes\psi_{ij}\\
w_1\otimes\psi_1\\
w_2\otimes\psi_2\\
\end{array}\right)
\mapsto\qt\left(\begin{array}{c}
2v_{220}\otimes(-\partial_1\psi_{21}-\partial_2\psi_{12})\\
2v_{202}\otimes(-\partial_1\psi_{31}-\partial_3\psi_{13})\\
2v_{022}\otimes(-\partial_3\psi_{23}-\partial_2\psi_{32})\\
v_{211}\otimes(-\partial_2\psi_{13}-\partial_3\psi_{12}-\partial_1\psi_1)\\
v_{121}\otimes(-\partial_1\psi_{23}-\partial_3\psi_{21}-\partial_2\psi_2)\\
v_{112}\otimes(\partial_2\psi_{31}+\partial_1\psi_{32}-\partial_3\psi_1-\partial_3\psi_2)\\
\end{array}\right).
\end{equation}
Let
\begin{eqnarray*}
(3,1,1):w_{311}&=&e_1\wedge e_3e_1\wedge e_2e_1-e_1\wedge e_2e_1\wedge e_3e_1\\
(1,3,1):w_{131}&=&e_2\wedge e_3e_1\wedge e_2e_2-e_1\wedge e_2e_2\wedge e_3e_2\\
(1,1,3):w_{113}&=&e_2\wedge e_3e_1\wedge e_3e_3-e_1\wedge e_3e_2\wedge e_3e_3\\
(2,2,1):w_{221}&=&e_2\wedge e_3e_1\wedge e_2e_1+e_1\wedge e_3e_1\wedge e_2e_2-e_1\wedge e_2e_2\wedge e_3e_1-e_1\wedge e_2e_1\wedge e_3e_2\\
(2,1,2):w_{212}&=&e_2\wedge e_3e_1\wedge e_3e_1+e_1\wedge e_3e_1\wedge e_2e_3-e_1\wedge e_3e_2\wedge e_3e_1-e_1\wedge e_2e_1\wedge e_3e_3\\
(1,2,2):w_{122}&=&e_2\wedge e_3e_1\wedge e_3e_2+e_2\wedge e_3e_1\wedge e_2e_3-e_1\wedge e_3e_2\wedge e_3e_2-e_1\wedge e_2e_2\wedge e_3e_3\\
\end{eqnarray*}
be a preferred basis of $\V_{\lambda_5}$. With respect to this basis
the operator $D_4$ is equal to
\begin{eqnarray*}
&\left(\begin{array}{c}
h_{ij}\otimes\psi_{ij}\\
w_1\otimes\psi_1\\
w_2\otimes\psi_2\\
\end{array}\right)
\mapsto&\\
&\frac{1}{8}\left(\begin{array}{c}
2w_{311}\otimes(\p_1\p_1(\psi_1+2\psi_2)+(2\p_2\p_1+\p_1\p_2)\psi_{13}-(\p_1\p_3+2\p_3\p_1)\psi_{12})\\
2w_{131}\otimes(-\p_2\p_2(2\psi_1+\psi_2)+(\p_2\p_3+2\p_3\p_2)\psi_{21}-(\p_2\p_1+2\p_1\p_2)\psi_{23})\\
2w_{113}\otimes(\p_3\p_3(\psi_1-\psi_2)+(\p_3\p_1+2\p_1\p_3)\psi_{32}-(\p_3\p_2+2\p_2\p_3)\psi_{31})\\
w_{221}\otimes(\p_1(\p_2\psi_2-3\p_1\psi_{23})+\p_2(-\p_1\psi_1+3\p_2\psi_{13})\\
+(\p_1\p_3+2\p_3\p_1)\psi_{21}-(\p_2\p_3+2\p_3\p_2)\psi_{12})\\
w_{212}\otimes(\p_1(\p_3\psi_1+\p_3\psi_2+3\p_1\psi_{32})+\p_3(\p_1\psi_1-3\p_3\psi_{12})\\
-(\p_1\p_2+2\p_2\p_1)\psi_{31}+(\p_3\p_2+2\p_2\p_3)\psi_{13})\\
w_{122}\otimes(\p_2(-\p_3\psi_1-\p_3\psi_2-3\p_2\psi_{31})+\p_3(-\p_2\psi_2+3\p_3\psi_{21})\\
+(\p_2\p_1+2\p_1\p_2)\psi_{32}-(\p_3\p_1+2\p_1\p_3)\psi_{23})
\end{array}\right).&
\end{eqnarray*}
Let us denote by $\{h_{ij}^\ast,w_1^\ast,w_2^\ast\}$ the basis of
$\V_{\lambda_6}$ dual to the basis of $\V_{\lambda_3}$. With respect
to this basis we can write the operator $D_5$ as
\begin{eqnarray*}
\left(\begin{array}{c}
w_{311}\otimes\psi_{311}\\
w_{131}\otimes\psi_{131}\\
w_{113}\otimes\psi_{113}\\
w_{221}\otimes\psi_{221}\\
w_{212}\otimes\psi_{212}\\
w_{122}\otimes\psi_{122}
\end{array}\right)&\mapsto&\frac{1}{24}
\left(\begin{array}{c}
h^\ast_{12}\otimes(\partial_3\psi_{122}-\partial_2\psi_{113})\\
h^\ast_{21}\otimes(\partial_3\psi_{212}-\partial_1\psi_{113})\\
h^\ast_{13}\otimes(\partial_2\psi_{122}-\partial_3\psi_{131})\\
h^\ast_{31}\otimes(\partial_2\psi_{221}-\partial_1\psi_{131})\\
h^\ast_{23}\otimes(\partial_1\psi_{212}-\partial_3\psi_{311})\\
h^\ast_{32}\otimes(\partial_1\psi_{221}-\partial_2\psi_{311})\\
w_{1}^\ast\otimes(\partial_1\psi_{122}-\partial_3\psi_{221})\\
w_{2}^\ast\otimes(\partial_2\psi_{212}-\partial_3\psi_{221})
\end{array}\right).
\end{eqnarray*}
The operator $D_6$ is given by
\begin{eqnarray*}
&\left(\begin{array}{c}
v_{220}\otimes\varphi_{220}\\
v_{202}\otimes\varphi_{202}\\
v_{022}\otimes\varphi_{022}\\
v_{211}\otimes\varphi_{211}\\
v_{121}\otimes\varphi_{121}\\
v_{112}\otimes\varphi_{112}
\end{array}\right)&\mapsto\\
&\frac{1}{48}
\left(\begin{array}{c}
h^\ast_{12}\otimes(-(\partial_3\partial_1+2\partial_1\partial_3)\varphi_{022}-(\partial_3\partial_2+2\partial_2\partial_3)\varphi_{112}+3\partial_3\partial_3\varphi_{121})\\
h^\ast_{21}\otimes((\partial_3\partial_2+2\partial_2\partial_3)\varphi_{202}+(\partial_3\partial_1+2\partial_1\partial_3)\varphi_{112}-3\partial_3\partial_3\varphi_{211})\\
h^\ast_{13}\otimes((2\partial_1\partial_2+\partial_2\partial_1)\varphi_{022}-(\partial_2\partial_3+2\partial_3\partial_2)\varphi_{121}+3\partial_2\partial_2\varphi_{112})\\
h^\ast_{31}\otimes(-(\partial_2\partial_3+2\partial_3\partial_2)\varphi_{220}+(\partial_2\partial_1+2\partial_1\partial_2)\varphi_{121}+3\partial_2\partial_2\varphi_{211})\\
h^\ast_{23}\otimes(-(\partial_1\partial_2+2\partial_2\partial_1)\varphi_{202}+(\partial_1\partial_3+2\partial_3\partial_1)\varphi_{211}-3\partial_1\partial_1\varphi_{112})\\
h^\ast_{32}\otimes((\partial_1\partial_3+2\partial_3\partial_1)\varphi_{220}-(\partial_1\partial_2+2\partial_2\partial_1)\varphi_{211}-3\partial_1\partial_1\varphi_{121})\\
w_{1}^\ast\otimes(\partial_1\partial_1\varphi_{022}-2\partial_2\partial_2\varphi_{202}+\partial_3\partial_3\varphi_{220}-\partial_1\partial_2\varphi_{112}-(\partial_3\partial_1+\partial_3\partial_1)\varphi_{121}+\partial_3\partial_2\varphi_{211})\\
w_{2}^\ast\otimes(2\partial_1\partial_1\varphi_{022}-\partial_2\partial_2\varphi_{202}-\partial_3\partial_3\varphi_{220}+\partial_2\partial_1\varphi_{112}-\partial_3\partial_1\varphi_{121}+(\partial_3\partial_2+\partial_2\partial_3)\varphi_{211})\\
\end{array}\right).&
\end{eqnarray*}
The operator $D_7$ is given by
\begin{eqnarray*}
\left(\begin{array}{c}
h^\ast_{12}\otimes\varphi_{12}\\
h^\ast_{21}\otimes\varphi_{21}\\
h^\ast_{13}\otimes\varphi_{13}\\
h^\ast_{31}\otimes\varphi_{31}\\
h^\ast_{23}\otimes\varphi_{23}\\
h^\ast_{32}\otimes\varphi_{32}\\
w_{1}^\ast\otimes\varphi_{1}\\
w_{2}^\ast\otimes\varphi_{2}\\
\end{array}\right)&\mapsto&
\left(\begin{array}{c}
e^\ast_1\otimes(\partial_2\partial_2\varphi_{21}-\partial_1\partial_2\varphi_{12}+\partial_3\partial_3\varphi_{31}-\partial_3\partial_1\varphi_{13}+\{\partial_2,\partial_3\}\varphi_1-\partial_3\partial_2\varphi_2)\\
e^\ast_2\otimes(\partial_2\partial_2\varphi_{12}-\partial_2\partial_1\varphi_{21}+\partial_3\partial_3\varphi_{32}-\partial_2\partial_3\varphi_{23}+\{\partial_1,\partial_3\}\varphi_2-\partial_3\partial_1\varphi_1)\\
e^\ast_3\otimes(\partial_1\partial_1\varphi_{13}-\partial_3\partial_1\varphi_{31}+\partial_2\partial_2\varphi_{23}-\partial_3\partial_2\varphi_{32}-\partial_2\partial_1\varphi_1-\partial_1\partial_2\varphi_2)\\
\end{array}\right).
\end{eqnarray*}
The operator $D_8$ is given by
\begin{eqnarray*}
\left(\begin{array}{c}
e^\ast_1\otimes\psi_1\\
e^\ast_2\otimes\psi_2\\
e^\ast_3\otimes\psi_3\\
\end{array}\right)\mapsto
\left(\begin{array}{c}
\partial_1\psi_1+\partial_2\psi_2+\partial_3\psi_3\\
\end{array}\right).
\end{eqnarray*}

\end{document}